\documentclass[11pt,reqno]{amsart}

\usepackage[T1]{fontenc}
\usepackage{lmodern}
\usepackage{microtype}
\usepackage{amsmath,amssymb,amsthm,mathtools}
\usepackage{booktabs}
\usepackage{hyperref}

\hypersetup{
  pdftitle={Bass-Serre Coverings, Profinite Completions, and Derivations of A\_l-Loops},
  pdfauthor={Raffaello Caserta, PhD},
  pdfsubject={Loop theory, Bass-Serre theory, profinite groups, fusion systems},
  pdfkeywords={A\_l-loop, left automorphic loop, left-envelope, fusion system, Bass-Serre theory, virtually free group, profinite completion, derivation},
  colorlinks=true,
  linkcolor=blue,
  citecolor=blue,
  urlcolor=blue
}

\newtheorem{theorem}{Theorem}[section]
\newtheorem{proposition}[theorem]{Proposition}
\newtheorem{lemma}[theorem]{Lemma}
\newtheorem{corollary}[theorem]{Corollary}

\theoremstyle{definition}
\newtheorem{definition}[theorem]{Definition}
\newtheorem{example}[theorem]{Example}
\newtheorem{problem}[theorem]{Open Problem}

\theoremstyle{remark}
\newtheorem{remark}[theorem]{Remark}

\DeclareMathOperator{\Aut}{Aut}
\DeclareMathOperator{\Sym}{Sym}
\DeclareMathOperator{\Hom}{Hom}
\DeclareMathOperator{\Inj}{Inj}
\DeclareMathOperator{\Mor}{Mor}
\DeclareMathOperator{\Stab}{Stab}
\DeclareMathOperator{\core}{core}
\DeclareMathOperator{\rk}{rk}
\DeclareMathOperator{\id}{id}

\newcommand{\F}{\mathcal{F}}
\newcommand{\G}{\mathcal{G}}
\newcommand{\T}{\mathcal{T}}
\newcommand{\Z}{\mathbb{Z}}
\newcommand{\inn}{\mathrm{inn}}

\title[Bass--Serre coverings and derivations of $A_\ell$-loops]{Bass--Serre coverings, profinite completions, and derivations of $A_\ell$-loops}

\author{Raffaello Caserta}
\address{Istituto Superiore ``Ettore Majorana'', Palermo, Italy}
\email{raffaello.caserta@majoranapa.edu.it}


\begin{document}

\begin{abstract}
We develop a group-theoretic and profinite framework connecting finite faithful left-envelopes $(G,H,L)$ of left-loops, fusion systems over $H$, and Bass--Serre theory. Given any finite faithful left-envelope $(G,H,L)$, we construct a discrete, virtually free, residually finite uniform tree lattice $\overline{G}=\pi_1(\G,\Gamma)$ acting faithfully on a $2(|L|-1)$-regular tree together with a left-envelope epimorphism $\psi\colon (\overline{G},H,\overline{L})\to(G,H,L)$ whose kernel $K=\ker\psi$ is a free normal subgroup of both $\overline{G}$ and $\overline{L}$ of exact rank $(|L|-1)^2$. We prove that $\psi$ preserves the loop property, the left inverse property (LIP), the $A_\ell$-property, and the fusion system $\F_H(\overline{G})=\F_H(G)$. Passing to the profinite completion $\widehat{\overline{G}}$ yields a continuous epimorphism onto $G$, a clopen section $\widehat{L}=\widehat{\psi}^{-1}(L)$ endowing $\widehat{L}$ with the structure of a profinite left-loop (and a profinite loop whenever $L$ is a loop), and a canonical projective tower of finite faithful left-envelopes covering $(G,H,L)$ with the same fusion system.

In the $A_\ell$-case, the fusion equality $\F_H(G)=\F_H(H)$ forces $\overline{G}$ to split as a semidirect product $\overline{G}=T\rtimes H$, where $T$ is a free group of rank $|L|-1$ on which $H$ acts by permuting a free basis. We establish an exact necessary and sufficient criterion for an inner derivation $F(T)=\{af(a)^{-1}:a\in T\}$ of an arbitrary semidirect product $T\rtimes H$ to be a loop. We then prove a central obstruction: any derivation factoring through $T/[T,H]$ that yields an $A_\ell$-loop necessarily has an abelian left inner mapping group $H(L_f)$, and if it is $H$-equivariant ($L_f^H=L_f$) then $f(T)\subseteq Z(H)$. Finally, passing to finite quotients of $T\rtimes H$ (in particular wreath products $(\Z/m\Z)^{|H|}\rtimes H$), we construct explicit infinite families of proper finite non-Bol $A_\ell$-loops with LIP and prescribed finite abelian left inner mapping group, exhibiting concrete examples of orders $9$, $16$, $25$, and $64$.
\end{abstract}

\maketitle

% ============================================================================
\section{Introduction}
% ============================================================================

Given a group $G$, a subgroup $H\le G$ (necessarily non-normal in the non-associative faithful case), and a left transversal $L$ to $H$ in $G$ containing the identity $1\in G$, it is elementary to see that $L$ inherits the structure of a left-loop $(L,\oplus)$. Conversely, any left-loop can be identified with a left transversal of the group $G(L)$ generated by its left translations $x\mapsto a\oplus x$ modulo the stabilizer $H(L)$ of the identity. Following a consolidated terminology\footnote{See, e.g., \cite[p.~17]{nagy2002loops}, \cite{aschbacher2005bol,aschbacher2006projective}, and \cite{sabinin1972equivalence}.}, we call a triple $(G,H,L)$ a \emph{left-folder}, and a \emph{faithful left-envelope} when $L$ generates $G$ and the core $\core_G(H)$ of the action of $G$ on $G/H$ is trivial. According to a classical result of Baer~\cite{baer1939nets} (cf.\ \cite[Proposition~1.6, p.~18]{nagy2002loops}), a left transversal $L$ is actually a loop precisely when it is simultaneously a left transversal to all conjugates $H^g$, $g\in G$.

Among non-associative loops, a prominent role is played by \emph{$A_\ell$-loops} (also termed \emph{left automorphic loops}), that is, loops whose left inner mappings $\delta_{a,b}=\lambda_{a\oplus b}^{-1}\lambda_a\lambda_b$ act as automorphisms of $(L,\oplus)$~\cite{bruck1956loops,kinyon2000loops,kinyon2016structure}. Notable instances are furnished by all \emph{Bruck loops} (left Bol loops satisfying the automorphic inverse property, or $K$-loops~\cite{kiechle2002theory,aschbacher2006bruck,baumeister2011bruck}), whereas general Bol loops---such as the finite Bol loops of exponent~$2$ studied in~\cite{baumeister2010self}---need not fulfill the $A_\ell$-condition. From a group-theoretical standpoint, the $A_\ell$-property of a faithful left-envelope $(G,H,L)$ is equivalent to the invariance $L^H=L$ of the transversal under conjugation by $H$~\cite[Lemma~4.5]{kinyon2000loops}, and it turns out that $H$ then controls $G$-fusion in $H$~\cite[(4.3)]{aschbacher2005bol}, i.e., $\F_H(G)=\F_H(H)$.

It is worth emphasizing that, given a finite group $G$, finding a non-normal subgroup $H<G$ admitting a left transversal $L$ which is a proper (non-associative) loop is a delicate problem, whereas constructing finite proper $A_\ell$-left-loops is comparatively frequent. The aim of this work is therefore to study finite faithful left-envelopes through the lens of Bass--Serre theory~\cite{serre1980trees,scott1979topological,bass1990uniform,bass2001tree} and profinite completions~\cite{ribes2010profinite}, showing how every finite left-envelope lifts canonically to both a discrete virtually free cover and a profinite tower, and applying the resulting semidirect structure in the $A_\ell$-case to build explicit finite and infinite proper $A_\ell$-loops with the left inverse property (LIP).

The paper is organized as follows. In Section~\ref{sec:preliminaries} we collect, for the convenience of the reader, the basic properties of faithful left-envelopes, fusion systems, and the $A_\ell$- and LIP conditions under uniform left-conjugation conventions, and establish a criterion (Proposition~\ref{prop:normal-quotient}) governing quotients of left-folders modulo normal subgroups saturating the transversal ($LK=L$). Section~\ref{sec:bass-serre} constructs the universal Bass--Serre covering $(\overline{G},H,\overline{L})$ of a finite faithful left-envelope $(G,H,L)$ (Theorem~\ref{thm:bass-serre-cover}), computes the rank of its free kernel $K$ via the rational Euler characteristic, and passes to the profinite completion $\widehat{\overline{G}}$ to obtain a profinite (left-)loop $\widehat{L}$ together with a projective tower of finite covers preserving $\F_H(G)$ (Theorem~\ref{thm:profinite-tower}); we then show that fusion control in the $A_\ell$-case forces a semidirect splitting $\overline{G}=T\rtimes H$ with $T$ free (Theorem~\ref{thm:al-semidirect}). Finally, Section~\ref{sec:derivations} analyzes inner derivations $F(T)=\{af(a)^{-1}:a\in T\}$ of an arbitrary semidirect product $T\rtimes H$: after proving an exact loop criterion (Theorem~\ref{thm:derivation-iff}), we identify a central obstruction for derivations factoring through the coinvariant quotient $T/[T,H]$ (Theorem~\ref{thm:obstruction}), and pass to finite quotients (Theorem~\ref{thm:finite-quotients} and Corollary~\ref{cor:wreath-family}) to produce explicit proper finite non-Bol $A_\ell$-loops with LIP (Example~\ref{ex:finite-wreath}), concluding with an open problem on non-factoring derivations for non-abelian $H$.

% ============================================================================
\section{Left-envelopes, normal quotients, and fusion systems}\label{sec:preliminaries}
% ============================================================================

Throughout the paper, all group conjugations act on the \emph{left}: for elements $x,g$ of a group $G$ and a subgroup $K\le G$, we write
\[
  x^g := g x g^{-1}, \qquad K^g := g K g^{-1}, \qquad \inn_g(x) := g x g^{-1}.
\]
Note that $(x^h)^g = x^{gh}$ and $\inn_g \circ \inn_h = \inn_{gh}$.

\subsection{Left-loops and faithful left-envelopes}

A \emph{left-loop} $(L,\oplus)$ is a set $L$ endowed with a binary operation $\oplus$ and a two-sided identity element $1\in L$ such that for all $a,b\in L$ the equation $a\oplus x = b$ has a unique solution $x\in L$. A \emph{right-loop} is defined dually via $x\oplus a = b$, and a \emph{loop} is both a left- and right-loop~\cite{bruck1958survey}. For each $a\in L$, the bijection $\lambda_a\colon x\mapsto a\oplus x$ is the \emph{left translation} by $a$. The subgroup $G(L):=\langle \lambda_a : a\in L\rangle \le \Sym(L)$ is the \emph{left multiplication} (or \emph{enveloping}) group of $L$, and the point stabilizer $H(L):=\Stab_{G(L)}(1)$ is its \emph{left inner mapping group}. Setting $T(L):=\{\lambda_a : a\in L\}$ and $\delta_{a,b}:=\lambda_{a\oplus b}^{-1}\lambda_a\lambda_b\in H(L)$, one has $G(L)=T(L)H(L)$ and $H(L)=\langle \delta_{a,b} : a,b\in L\rangle$~\cite[Section~1.2]{nagy2002loops}.

A \emph{left-folder} is a triple $(G,H,L)$ consisting of a group $G$, a subgroup $H\le G$, and a left transversal $L$ to $H$ in $G$ containing $1\in G$ (equivalently, $L=\sigma(G/H)$ for a \emph{section} $\sigma\colon G/H\to G$ satisfying $\pi\circ\sigma=\id_{G/H}$ and $\sigma(H)=1$, where $\pi\colon G\to G/H$ is the canonical projection; see~\cite[p.~17]{nagy2002loops}). The transversal $L$ inherits a left-loop operation $a\oplus b := c$, where $\{c\} = abH\cap L$. A left-folder is a \emph{left-envelope} if $G=\langle L\rangle$, and is \emph{faithful} if the core $\core_G(H):=\bigcap_{g\in G} H^g$ of the left action of $G$ on $G/H$ is trivial.

A \emph{morphism} of left-folders $\varphi\colon (G_1,H_1,L_1)\to (G_2,H_2,L_2)$ is a group homomorphism $\varphi\colon G_1\to G_2$ such that $\varphi(H_1)\le H_2$ and $\varphi(L_1)\subseteq L_2$; it is an \emph{epimorphism} if $\varphi(G_1)=G_2$, $\varphi(H_1)=H_2$, and $\varphi(L_1)=L_2$. For any left-folder $(G,H,L)$, the left action of $G$ on $G/H\cong L$ defines a homomorphism $\lambda\colon G\to\Sym(L)$ via $\{\lambda(g)(x)\} := gxH\cap L$ ($g\in G$, $x\in L$), with kernel $\core_G(H)$. By~\cite[Theorem~1.11, pp.~21--22]{nagy2002loops} (see also~\cite{aschbacher2005bol,foguel2000groups,lal1996transversals,niemenmaa1990multiplication}), if $(G,H,L)$ is a faithful left-envelope, then $\lambda$ restricts to a left-envelope isomorphism $(G,H,L)\xrightarrow{\sim} (G(L),H(L),T(L))$.

Following Kinyon and Jones~\cite[Section~2]{kinyon2000loops}, for a left-folder $(G,H,L)$ every element of $G$ writes uniquely as $ah$ ($a\in L$, $h\in H$). For $a,b\in L$ and $h\in H$, define
\[
  (a,b) := (a\oplus b)^{-1}ab \in H, \quad
  \{h\cdot a\} := haH\cap L, \quad
  \mu_a(h) := (h\cdot a)^{-1}ha \in H.
\]
Then $ab = (a\oplus b)(a,b)$ and $ha = (h\cdot a)\mu_a(h)$. Notice that $\lambda((a,b))=\delta_{a,b}$ and $\lambda(h)(a)=h\cdot a$, and since $H$ normalizes $H$ we have $\{h\cdot a\} = a^h H\cap L$, whence the stabilizer of $a\in L$ in $H$ is
\[
  H_a := \Stab_H(a) = \{h\in H : haH = aH\} = H\cap H^a.
\]

\begin{proposition}\label{prop:basic-identities}
Let $(G,H,L)$ be a left-folder. For all $a,b\in L$ and $h,k\in H$, one has:
\begin{enumerate}
  \item $h\cdot 1 = 1$ and $(hk)\cdot a = h\cdot(k\cdot a)$;
  \item $h\cdot(a\oplus b) = (h\cdot a)\oplus(\mu_a(h)\cdot b)$.
\end{enumerate}
\end{proposition}

\begin{proof}
Part~(1) is immediate from $h1H = H = 1H$ and $(hk)aH = h(k\cdot a)H = (h\cdot(k\cdot a))H$. For~(2), comparing the unique $L$-components of $hab = (h\cdot(a\oplus b))\,\mu_{a\oplus b}(h)\,(a,b)$ and $hab = (h\cdot a)(\mu_a(h)\cdot b)\,\mu_b(\mu_a(h)) = ((h\cdot a)\oplus(\mu_a(h)\cdot b))\,(h\cdot a,\mu_a(h)\cdot b)\,\mu_b(\mu_a(h))$ in $G=LH$ yields $h\cdot(a\oplus b) = (h\cdot a)\oplus(\mu_a(h)\cdot b)$.
\end{proof}

\begin{definition}\label{def:al-lip}
Let $(L,\oplus)$ be a left-loop.
\begin{enumerate}
  \item $(L,\oplus)$ is an \emph{$A_\ell$-left-loop} (or \emph{left automorphic left-loop}) if $H(L)\le \Aut(L,\oplus)$. A faithful left-envelope $(G,H,L)$ whose underlying left-loop $(L,\oplus)$ is an $A_\ell$-left-loop is called a \emph{faithful $A_\ell$-left-envelope}.
  \item $(L,\oplus)$ satisfies the \emph{left inverse property} (\emph{LIP}) if for every $a\in L$ there exists $a^\lambda\in L$ such that $a^\lambda\oplus(a\oplus b) = b$ for all $b\in L$.
\end{enumerate}
\end{definition}

The group-theoretic characterizations of LIP and of the $A_\ell$-property are due to Kinyon and Jones~\cite[Corollary~4.3 and Lemma~4.5]{kinyon2000loops}; we record the faithful case for convenient reference.

\begin{proposition}\label{prop:al-lip-char}
Let $(G,H,L)$ be a faithful left-envelope.
\begin{enumerate}
  \item $(L,\oplus)$ satisfies the LIP if and only if $L^{-1} = L$ (and then $a^\lambda = a^{-1}$ in $G$ for all $a\in L$).
  \item The following are equivalent:
  \begin{enumerate}
    \item[(i)] $(L,\oplus)$ is an $A_\ell$-left-loop;
    \item[(ii)] $L^h = L$ for all $h\in H$ (and then $h\cdot a = a^h$);
    \item[(iii)] $\mu_a = \id_H$ for all $a\in L$.
  \end{enumerate}
\end{enumerate}
\end{proposition}

\begin{proof}
Part~(1) is~\cite[Corollary~4.3]{kinyon2000loops}: $(L,\oplus)$ has LIP if and only if $\lambda_a^{-1}\in T(L)$ for all $a\in L$, which under the isomorphism $G\cong G(L)$ is precisely $L^{-1}=L$.

For~(2), if $H(L)\le \Aut(L,\oplus)$, then for all $h\in H$ and $a,x\in L$, $\lambda(h)(a\oplus x) = \lambda(h)(a)\oplus \lambda(h)(x)$, i.e., $\lambda(h)\lambda_a\lambda(h)^{-1} = \lambda_{\lambda(h)(a)}$. Pulling back via the isomorphism $G\cong G(L)$ gives $a^h = h\cdot a \in L$, so~(i)$\Rightarrow$(ii). If~(ii) holds, then $a^h\in L\cap haH = \{h\cdot a\}$, whence $\mu_a(h) = (a^h)^{-1}ha = h$, giving~(ii)$\Rightarrow$(iii). Finally, if $\mu_a=\id_H$ for all $a\in L$, Proposition~\ref{prop:basic-identities}(2) gives $h\cdot(a\oplus b) = (h\cdot a)\oplus(h\cdot b)$ for all $h\in H$, so $H(L)\le\Aut(L,\oplus)$, giving~(iii)$\Rightarrow$(i).
\end{proof}

\subsection{Fusion systems and \texorpdfstring{$A_\ell$}{A\_l}-envelopes}

Following~\cite[Definition~1.1]{broto2003homotopy} (dropping the finite $p$-group restriction), we use the categorical language of fusion systems over an arbitrary group $H$.

\begin{definition}\label{def:fusion-system}
A \emph{fusion system} $\F$ over a group $H$ is a category whose objects are the subgroups of $H$ and whose morphism sets $\Mor_\F(R,S)$ satisfy:
\begin{enumerate}
  \item $\Hom_H(R,S) \subseteq \Mor_\F(R,S) \subseteq \Inj(R,S)$, where $\Hom_H(R,S) = \{\inn_h|_R : h\in H,\; R^h\le S\}$;
  \item every $\phi\in\Mor_\F(R,S)$ factors as an isomorphism $\phi\colon R\xrightarrow{\sim}\phi(R)$ in $\F$ followed by the inclusion $\phi(R)\hookrightarrow S$, and $\phi^{-1}\in\Mor_\F(\phi(R),R)$.
\end{enumerate}
If $G$ is a group containing $H$, the \emph{fusion system of $G$ over $H$}, denoted $\F_H(G)$, has morphisms $\Mor_{\F_H(G)}(R,S) := \{\inn_g|_R : g\in G,\; R^g\le S\}$.
\end{definition}

\begin{proposition}[cf.\ {\cite[(4.3)]{aschbacher2005bol}}]\label{prop:al-fusion}
Let $(G,H,L)$ be a faithful left-envelope. If $(L,\oplus)$ is an $A_\ell$-left-loop, then for every $g=ah\in G$ ($a\in L$, $h\in H$) and every $t\in H\cap H^{g^{-1}}$, one has $t^g = t^h$. In particular, $\F_H(G) = \F_H(H)$, i.e., $H$ controls $G$-fusion in $H$.
\end{proposition}

\begin{proof}
Let $t\in H$ be such that $t' := t^g \in H$, where $g = ah$ with $a\in L$ and $h\in H$. Then $a t^h a^{-1} = t'$, so
\[
  a\,t^h = t' a = (t' a t'^{-1})\,t' = a^{t'}\,t'.
\]
By Proposition~\ref{prop:al-lip-char}(2), $a^{t'}\in L$. By uniqueness of the decomposition $G = LH$, we deduce $a^{t'} = a$ and $t^h = t' = t^g$. Thus $\inn_g|_R = \inn_h|_R$ whenever $R\le H$ and $R^g\le H$, proving $\F_H(G)\subseteq \F_H(H)\subseteq \F_H(G)$.
\end{proof}

\subsection{Normal subgroups contained in a transversal}

A recurring step in our covering and quotient constructions is passing modulo a normal subgroup $K\trianglelefteq G$ that saturates the transversal $L$. The following proposition summarizes how the structural properties of a left-folder behave under such quotients.

\begin{proposition}\label{prop:normal-quotient}
Let $(G,H,L)$ be a left-folder and let $K\trianglelefteq G$ be a normal subgroup of $G$ such that $LK = L$ (equivalently, $K\subseteq L$ and $L$ is a union of cosets of $K$). Let $\pi_K\colon G\to G/K$ be the canonical projection, and set $\widetilde{G}:=G/K$, $\widetilde{H}:=HK/K\cong H$, and $\widetilde{L}:=L/K = \{aK : a\in L\}$. Then:
\begin{enumerate}
  \item $K\cap H = \{1\}$, and $(a,k) = (k,a) = 1$ for all $a\in L$ and $k\in K$. In particular, $k_1\oplus k_2 = k_1k_2$ on $K$, so $K$ is a subgroup of the left-loop $(L,\oplus)$ (i.e., a sub-left-loop forming a group under $\oplus$).
  \item $(a\oplus k, b) = (a,b) = (a, b\oplus k)$ for all $a,b\in L$ and $k\in K$.
  \item $(\widetilde{G},\widetilde{H},\widetilde{L})$ is a left-folder (and a left-envelope whenever $(G,H,L)$ is), and $\pi_K|_L\colon (L,\oplus)\to(\widetilde{L},\oplus)$ is a left-loop epimorphism with kernel $K$ (whence $K$ is a normal sub-left-loop, indeed a normal subgroup, of $(L,\oplus)$).
  \item $(L,\oplus)$ is a loop if and only if $(\widetilde{L},\oplus)$ is a loop.
  \item $L^H = L$ if and only if $\widetilde{L}^{\widetilde{H}} = \widetilde{L}$. In particular, if both $(G,H,L)$ and $(\widetilde{G},\widetilde{H},\widetilde{L})$ are faithful, then $(L,\oplus)$ is an $A_\ell$-left-loop if and only if $(\widetilde{L},\oplus)$ is.
  \item $L^{-1} = L$ if and only if $\widetilde{L}^{-1} = \widetilde{L}$. In particular, if both envelopes are faithful, $(L,\oplus)$ satisfies the LIP if and only if $(\widetilde{L},\oplus)$ does.
\end{enumerate}
\end{proposition}

\begin{proof}
(1) Since $K\subseteq L$ and $L$ is a left transversal to $H$ containing $1$, $K\cap H\subseteq L\cap H = \{1\}$. For $a\in L$ and $k\in K$, we have $ak\in LK = L$, hence $a\oplus k = ak$ and $(a,k) = 1$. Since $K\trianglelefteq G$, $ka = a(a^{-1}ka)\in LK = L$, hence $k\oplus a = ka$ and $(k,a)=1$.

(2) For $a,b\in L$ and $k\in K$, using normality of $K$ we write
\[
  (a\oplus k)b = akb = ab(b^{-1}kb) = (a\oplus b)(a,b)(b^{-1}kb) = (a\oplus b)\,k'\,(a,b),
\]
where $k':=(a,b)(b^{-1}kb)(a,b)^{-1}\in K$. Since $(a\oplus b)k'\in LK = L$, uniqueness of the decomposition $G = LH$ implies $(a\oplus k)\oplus b = (a\oplus b)k'$ and $(a\oplus k,b) = (a,b)$. Similarly,
\[
  a(b\oplus k) = abk = (a\oplus b)(a,b)k = (a\oplus b)\,k''\,(a,b)
\]
with $k'' = (a,b)k(a,b)^{-1}\in K$, yielding $(a,b\oplus k) = (a,b)$.

(3) Since $G = LH$, we have $\widetilde{G} = \widetilde{L}\widetilde{H}$. If $aK,bK\in\widetilde{L}$ satisfy $aK\widetilde{H} = bK\widetilde{H}$, then $a\in bKH = bHK$, so $a = bhk = (bk')h$ for some $h\in H$ and $k,k'\in K$. Since $bk'\in LK=L$ and $L$ is a transversal to $H$, $h=1$ and $a=bk'\in bK$, so $aK = bK$. Thus $\widetilde{L}$ is a left transversal to $\widetilde{H}$ in $\widetilde{G}$ (and if $\langle L\rangle = G$, then $\langle\widetilde{L}\rangle = \pi_K(\langle L\rangle) = \widetilde{G}$). Clearly $\pi_K(a\oplus b) = \pi_K(ab(a,b)^{-1}) = \pi_K(a)\pi_K(b)\pi_K(a,b)^{-1}\in \widetilde{L}$ implies $\pi_K(a\oplus b) = \pi_K(a)\oplus \pi_K(b)$ and $\ker(\pi_K|_L) = L\cap K = K$.

(4) Recall that by~\cite[Proposition~1.6]{nagy2002loops} (or directly by solving $x\oplus a = b$), $(L,\oplus)$ is a loop if and only if for every $a\in L$ the right translation $\rho_a\colon x\mapsto x\oplus a$ is bijective on $L$. Note first that for all $x,a\in L$ and $k\in K$, by~(1) and~(2) we have
\begin{equation}\label{eq:coset-right-trans}
  (x\oplus k)\oplus a = (x\oplus a)\oplus k^*, \qquad k^* := (x,a)(a^{-1}ka)(x,a)^{-1} \in K.
\end{equation}
For fixed $x,a\in L$, the map $k\mapsto k^*$ is an automorphism of the group $K$ (conjugation by $(x,a)a^{-1}\in G$).

Assume $(\widetilde{L},\oplus)$ is a loop, and let $a,b\in L$. For surjectivity of $\rho_a$, choose $x_0K\in\widetilde{L}$ such that $(x_0K)\oplus(aK) = bK$. Then $x_0\oplus a = b\oplus k_0$ for some $k_0\in K$. Since $k\mapsto k^*$ is an automorphism of $K$, there exists $k\in K$ such that $k^* = k_0^{-1}$. By~\eqref{eq:coset-right-trans} and~(1),
\[
  (x_0\oplus k)\oplus a = (x_0\oplus a)\oplus k_0^{-1} = (b\oplus k_0)\oplus k_0^{-1} = b k_0 k_0^{-1} = b.
\]
For injectivity of $\rho_a$, suppose $x_1\oplus a = x_2\oplus a$. Projecting to $\widetilde{L}$ gives $(x_1K)\oplus(aK) = (x_2K)\oplus(aK)$, hence $x_1K = x_2K$ because $\widetilde{L}$ is a loop, so $x_1 = x_2\oplus k$ for some $k\in K$. By~\eqref{eq:coset-right-trans}, $x_2\oplus a = x_1\oplus a = (x_2\oplus a)\oplus k^* = (x_2\oplus a)k^*$, forcing $k^* = 1$, hence $k=1$ and $x_1=x_2$. Thus $L$ is a loop.

Conversely, assume $(L,\oplus)$ is a loop. Then $\widetilde{\rho}_{aK}\colon xK\mapsto (x\oplus a)K$ is manifestly surjective on $\widetilde{L}$. If $(x_1\oplus a)K = (x_2\oplus a)K$, then $x_1\oplus a = (x_2\oplus a)\oplus k_1$ for some $k_1\in K$. Choosing $k\in K$ with $k^* = k_1$ in~\eqref{eq:coset-right-trans} yields $(x_2\oplus k)\oplus a = (x_2\oplus a)\oplus k_1 = x_1\oplus a$, so $x_1 = x_2\oplus k \in x_2K$ by injectivity of $\rho_a$. Hence $\widetilde{L}$ is a loop.

(5) For $a\in L$ and $h\in H$, $(aK)^{hK} = a^h K\in \widetilde{L} \iff a^h\in LK = L$. Thus $\widetilde{L}^{\widetilde{H}} = \widetilde{L} \iff L^H = L$. The faithful statement follows from Proposition~\ref{prop:al-lip-char}(2).

(6) For $a\in L$, $(aK)^{-1} = a^{-1}K\in \widetilde{L} \iff a^{-1}\in LK = L$. The faithful statement follows from Proposition~\ref{prop:al-lip-char}(1).
\end{proof}

% ============================================================================
\section{Bass--Serre coverings and profinite completions}\label{sec:bass-serre}
% ============================================================================

\subsection{Flower graphs of groups and HNN normal forms}

Let $H$ be a finite group, let $S$ be a finite index set, and for each $s\in S$ let $H_s, H_s'\le H$ be subgroups equipped with an isomorphism $\partial_s\colon H_s\xrightarrow{\sim} H_s'$. Consider the \emph{flower graph of groups} $(\G,\Gamma)$ having a single vertex $v$ labeled by $\Gamma_v = H$ and $|S|$ oriented loops (edges) $\{e_s : s\in S\}$ based at $v$, with edge group $\Gamma_{e_s} = H_s$ and boundary monomorphisms $\partial_{e_s}^+ = \id_{H_s}\colon H_s\hookrightarrow H$ and $\partial_{e_s}^- = \partial_s\colon H_s\hookrightarrow H$.

The fundamental group $\overline{G} := \pi_1(\G,\Gamma)$ is the multiple HNN extension~\cite{higman1949embedding}
\begin{equation}\label{eq:hnn-presentation}
  \overline{G} = \bigl\langle H,\; \{e_s : s\in S\} \;\big|\; e_s h e_s^{-1} = \partial_s(h)\;\; \forall s\in S,\; \forall h\in H_s \bigr\rangle.
\end{equation}
Set $E^+ := \{e_s : s\in S\}$, $E^- := \{e_s^{-1} : s\in S\}$, and $E := E^+\sqcup E^-$. For $e = e_s^{-1}\in E^-$, set $H_e := H_s'$, $H_e' := H_s$, and $\partial_e := \partial_s^{-1}\colon H_s'\xrightarrow{\sim} H_s$, so that $e h e^{-1} = \partial_e(h)$ for all $e\in E$ and $h\in H_e$.

By Bass--Serre theory~\cite[Section~I.5.3]{serre1980trees} (see also~\cite[Section~4]{scott1979topological}, \cite[Section~1]{bass1990uniform}), $\overline{G}$ acts without inversions on its Bass--Serre tree $\T$, whose vertex set is $V(\T) = \overline{G}/H$ and where two vertices $g_1H, g_2H$ are adjacent if and only if $g_2H = g_1heH$ for some $h\in H$ and $e\in E$. The stabilizer of a vertex $gH$ is $H^g$, and the stabilizer of the directed edge from $gH$ to $gheH$ is $H^g\cap H^{ghe} = (H_e')^{gh}$. Since $\overline{G}$ acts vertex-transitively on $\T$, $\T$ is a $d$-regular locally finite tree of even degree
\begin{equation}\label{eq:tree-degree}
  d = \sum_{s\in S} \bigl([H:H_s] + [H:H_s']\bigr) = 2\sum_{s\in S}[H:H_s].
\end{equation}
The kernel of the action of $\overline{G}$ on $\T$ is $\core_{\overline{G}}(H) = \bigcap_{g\in\overline{G}} H^g$. Because the vertex stabilizers are finite and the quotient graph $\overline{G}\backslash\T = \G$ is finite, the faithful quotient $\overline{G}/\core_{\overline{G}}(H)$ is a discrete \emph{uniform tree lattice} in $\Aut(\T)$~\cite[Section~1]{bass1990uniform}, \cite[Chapter~1]{bass2001tree} (and $\overline{G}$ itself embeds whenever $\core_{\overline{G}}(H)=\{1\}$).

We recall two standard consequences of the HNN and Bass--Serre structure:
\begin{itemize}
  \item \textbf{Normal form (Britton's Lemma, \cite[Theorem~IV.2.1]{lyndon2001combinatorial}).} For each $e\in E$, fix a left transversal $R_e$ to $H_e'$ in $H$ containing $1$. Every element $g\in\overline{G}$ can be written \emph{uniquely} in the reduced form
  \begin{equation}\label{eq:hnn-normal-form}
    g = r_1 e_1 r_2 e_2 \cdots r_n e_n h, \qquad n\ge 0,\; e_i\in E,\; r_i\in R_{e_i},\; h\in H,
  \end{equation}
  subject to the non-backtracking condition that $r_{i+1}\ne 1$ whenever $e_{i+1} = e_i^{-1}$ ($1\le i < n$). (Indeed, $h_i e_i = r_i h_i' e_i = r_i e_i \partial_{e_i}^{-1}(h_i')$ pushes elements of $H_{e_i}'$ to the right across $e_i$.)
  \item \textbf{Virtual freeness and Euler characteristic.} A subgroup $K\le\overline{G}$ is free if and only if $K\cap H^g = \{1\}$ for all $g\in\overline{G}$~\cite[Section~I.3.3, Theorem~4]{serre1980trees}. Moreover, $\overline{G}$ possesses a normal free subgroup of finite index (i.e., $\overline{G}$ is \emph{virtually free}, hence \emph{residually finite}~\cite[Section~II.2.6, Proposition~11]{serre1980trees}), and its rational Euler characteristic~\cite[Section~7]{scott1979topological}, \cite[Section~1]{bass1990uniform} is
  \begin{equation}\label{eq:euler-char}
    \chi(\overline{G}) = \frac{1}{|H|} - \sum_{s\in S}\frac{1}{|H_s|}.
  \end{equation}
\end{itemize}

\subsection{The Bass--Serre covering of a finite left-envelope}

Let $(G,H,L)$ be a finite faithful left-envelope. The subgroup $H$ acts on $L^\# := L\setminus\{1\}$ via $h\cdot a \in haH\cap L$. Choose a complete set $S\subseteq L^\#$ of representatives of the $H$-orbits on $L^\#$. For each $s\in S$, set
\begin{align*}
  H_s &:= H\cap H^{s^{-1}}, \qquad H_s' := H\cap H^s = \Stab_H(s), \\
  \partial_s &:= \inn_s|_{H_s}\colon H_s\xrightarrow{\sim} H_s'.
\end{align*}
Notice that $|H_s| = |H_s'|$ and $[H:H_s] = [H:H_s'] = |\mathrm{Orb}_H(s)|$, whence
\begin{equation}\label{eq:orbit-sum}
  \sum_{s\in S} [H:H_s] = \sum_{s\in S} |\mathrm{Orb}_H(s)| = |L^\#| = |L| - 1.
\end{equation}
Let $\overline{G} := \pi_1(\G,\Gamma)$ be the fundamental group~\eqref{eq:hnn-presentation} attached to $(H,S,\{H_s,\partial_s\}_{s\in S})$.

\begin{theorem}[Bass--Serre universal covering]\label{thm:bass-serre-cover}
Let $(G,H,L)$ be a finite faithful left-envelope, and let $\overline{G}$ be constructed as above. Then:
\begin{enumerate}
  \item There exists a group epimorphism $\psi\colon \overline{G}\to G$ satisfying $\psi|_H = \id_H$ and $\psi(e_s) = s$ for all $s\in S$.
  \item Setting $\overline{L} := \psi^{-1}(L)\subseteq \overline{G}$, the triple $(\overline{G},H,\overline{L})$ is a faithful left-envelope, $\psi\colon (\overline{G},H,\overline{L})\to(G,H,L)$ is a left-envelope epimorphism, and $K := \ker\psi$ satisfies $K\subseteq \overline{L}$ and $\overline{L}K = \overline{L}$.
  \item The kernel $K$ is a free normal subgroup of both $\overline{G}$ and $(\overline{L},\oplus)$ of exact rank
  \[
    \rk(K) = (|L|-1)^2.
  \]
  \item The Bass--Serre tree $\T$ of $\overline{G}$ is $2(|L|-1)$-regular, and $\overline{G}$ is a discrete, virtually free, residually finite uniform lattice in $\Aut(\T)$.
  \item $(\overline{L},\oplus)$ is a loop (resp.\ an $A_\ell$-left-loop, resp.\ satisfies the LIP) if and only if $(L,\oplus)$ is.
  \item The fusion systems coincide: $\F_H(\overline{G}) = \F_H(G)$.
\end{enumerate}
\end{theorem}

\begin{proof}
(1) For each $s\in S$ and $h\in H_s = H\cap H^{s^{-1}}$, we have $shs^{-1} = \partial_s(h)$ in $G$. By the universal property of the presentation~\eqref{eq:hnn-presentation}, the assignments $h\mapsto h$ ($h\in H$) and $e_s\mapsto s$ ($s\in S$) extend uniquely to a group homomorphism $\psi\colon \overline{G}\to G$. To prove surjectivity, take any $a\in L^\#$. By definition of $S$, there exist $s\in S$ and $h\in H$ such that $a = h\cdot s$, so $a = hsh'$ for some $h'\in H$. Hence $a\in\langle H,S\rangle = \psi(\overline{G})$. Since $L$ generates $G$, we conclude $\psi(\overline{G}) = G$.

(2) Let $g\in\overline{G}$. Since $G=LH$, there exist unique $a\in L$ and $h\in H$ such that $\psi(g) = ah$. Then $\psi(gh^{-1}) = \psi(g)h^{-1} = a\in L$, so $gh^{-1}\in\overline{L}$ and $g = (gh^{-1})h\in\overline{L}H$. If $x,y\in\overline{L}$ satisfy $xH = yH$, then $x = yh$ for some $h\in H$, whence $\psi(x) = \psi(y)h$ with $\psi(x),\psi(y)\in L$; since $L$ is a transversal to $H$ in $G$, $h=1$ and $x=y$. Thus $\overline{L}$ is a left transversal to $H$ in $\overline{G}$ containing $1$ and $\{e_s : s\in S\}$, and clearly $\overline{L}K = \psi^{-1}(L) = \overline{L}$ with $K = \psi^{-1}(1)\subseteq\overline{L}$.

To verify $\langle\overline{L}\rangle = \overline{G}$, let $h\in H$. Since $L$ generates $G$, we can write $h = a_1\cdots a_r$ in $G$ with $a_i\in L$. Choosing preimages $\bar{a}_i\in\overline{L}$ with $\psi(\bar{a}_i)=a_i$, we have $k := (\bar{a}_1\cdots\bar{a}_r)^{-1}h \in \ker\psi = K \subseteq \overline{L}$, so $h = \bar{a}_1\cdots\bar{a}_r k\in\langle\overline{L}\rangle$. Thus $H\le\langle\overline{L}\rangle$, and hence $\overline{G} = \overline{L}H = \langle\overline{L}\rangle$. For faithfulness, if $h\in \core_{\overline{G}}(H) = \bigcap_{g\in\overline{G}} H^g$, then applying the epimorphism $\psi$ gives $h = \psi(h) \in \bigcap_{g\in\overline{G}} H^{\psi(g)} = \core_G(H) = \{1\}$.

(3) Since $\psi|_H = \id_H$, $K\cap H = \{1\}$. By normality of $K$, $K\cap H^g = (K\cap H)^g = \{1\}$ for every $g\in\overline{G}$, so $K$ acts freely on the Bass--Serre tree $\T$ and is therefore a free group. Because $\overline{G}$ is finitely generated and $[\overline{G}:K] = |G| < \infty$, $K$ is finitely generated; multiplicativity of the rational Euler characteristic~\cite[Section~7]{scott1979topological} together with~\eqref{eq:euler-char} and~\eqref{eq:orbit-sum} yields
\begin{align*}
  1 - \rk(K) = \chi(K) &= |G|\,\chi(\overline{G})
  = \frac{|G|}{|H|}\left(1 - \sum_{s\in S}[H:H_s]\right) \\
  &= |L|\bigl(1 - (|L|-1)\bigr) = |L|(2-|L|).
\end{align*}
Consequently, $\rk(K) = 1 + |L|(|L|-2) = (|L|-1)^2$. By Proposition~\ref{prop:normal-quotient}(1) and~(3), $K$ is also a normal subgroup of $(\overline{L},\oplus)$.

(4) Follows immediately from~\eqref{eq:tree-degree}, \eqref{eq:orbit-sum}, and the existence of the finite-index normal free subgroup $K\trianglelefteq \overline{G}$. Notice moreover that the kernel of the action of $\overline{G}$ on $\T$ is $\bigcap_{g\in\overline{G}} H^g = \core_{\overline{G}}(H) = \{1\}$ by~(2), so $\overline{G}$ acts faithfully on $\T$ and hence embeds as a discrete uniform lattice in $\Aut(\T)$.

(5) Since both $(\overline{G},H,\overline{L})$ and $(G,H,L)\cong (\overline{G}/K,HK/K,\overline{L}/K)$ are faithful left-envelopes and $\overline{L}K = \overline{L}$, part~(5) is a direct consequence of Proposition~\ref{prop:normal-quotient}(4)--(6).

(6) Since $\psi|_H = \id_H$, any conjugation by $g\in\overline{G}$ mapping a subgroup $R\le H$ into $H$ induces the same morphism as conjugation by $\psi(g)\in G$, so $\F_H(\overline{G})\subseteq \F_H(G)$. Conversely, let $g\in G$ and $R\le H$ satisfy $R^g\le H$. Write $g = ah_0$ with $a\in L$ and $h_0\in H$. If $a=1$, then $\inn_g|_R = \inn_{h_0}|_R\in\F_H(H)\subseteq\F_H(\overline{G})$. If $a\in L^\#$, write $a = h_1 s h_2$ with $s\in S$ and $h_1,h_2\in H$; setting $h_3 := h_2 h_0\in H$, we have $g = h_1 s h_3$. Then $(R^{h_3})^s = (R^g)^{h_1^{-1}}\le H$, so $R^{h_3}\le H\cap H^{s^{-1}} = H_s$. Lifting $g$ to $\bar{g} := h_1 e_s h_3 \in \overline{G}$ gives, for every $r\in R$,
\[
  r^{\bar{g}} = \bigl((r^{h_3})^{e_s}\bigr)^{h_1} = \bigl(\partial_s(r^{h_3})\bigr)^{h_1} = \bigl((r^{h_3})^s\bigr)^{h_1} = r^g.
\]
Hence $\inn_g|_R = \inn_{\bar{g}}|_R\in\Mor_{\F_H(\overline{G})}(R,H)$, proving $\F_H(G) = \F_H(\overline{G})$.
\end{proof}

\begin{remark}\label{rmk:not-profinite-not-lie}
It is important to distinguish the discrete group $\overline{G}$ from its topological completions. Although $\overline{G}$ embeds as a uniform lattice in the locally compact totally disconnected group $\Aut(\T)$, $\overline{G}$ itself is countably infinite and discrete for $|L|\ge 2$, hence non-compact and not profinite. Furthermore, for $d\ge 3$, $\Aut(\T)$ is not a $p$-adic Lie group in the sense of Lazard (its vertex stabilizers are iterated infinite wreath products), and the even degree $d=2(|L|-1)$ need not be of the form $p^m+1$. The genuine profinite object attached to $(G,H,L)$ is the \emph{profinite completion} $\widehat{\overline{G}}$, to which we now turn.
\end{remark}

\subsection{Profinite completions and projective towers of finite loops}

Let $\widehat{\overline{G}} := \varprojlim_{N\trianglelefteq_f \overline{G}} \overline{G}/N$ be the profinite completion of $\overline{G}$~\cite[Chapter~3]{ribes2010profinite}, where $N$ ranges over the directed system of finite-index normal subgroups of $\overline{G}$. We define a \emph{profinite (left-)loop} to be a topological (left-)loop (with continuous binary operation and left division, and continuous right division in the loop case) that is isomorphic as a topological (left-)loop to an inverse limit of finite discrete (left-)loops; in particular, every profinite (left-)loop is compact, Hausdorff, and totally disconnected\footnote{Whether every compact, Hausdorff, totally disconnected topological loop (or \emph{Boolean loop}) is necessarily profinite remains an open problem; see~\cite{herfort2011boolean}.}.

\begin{theorem}[Profinite covering and projective tower]\label{thm:profinite-tower}
Under the hypotheses of Theorem~\ref{thm:bass-serre-cover}, let $\widehat{\overline{G}}$ be the profinite completion of $\overline{G}$. Then:
\begin{enumerate}
  \item The canonical homomorphism $\iota\colon \overline{G}\to\widehat{\overline{G}}$ is injective with dense image, and $\iota(H)$ is a finite (hence closed) subgroup of $\widehat{\overline{G}}$ identified with $H$.
  \item The epimorphism $\psi\colon \overline{G}\to G$ extends uniquely to a continuous epimorphism $\widehat{\psi}\colon \widehat{\overline{G}}\to G$ whose kernel $\widehat{K} := \ker\widehat{\psi}$ is the closure of $K$ in $\widehat{\overline{G}}$, isomorphic to the free profinite group $\widehat{F}_{(|L|-1)^2}$ of rank $(|L|-1)^2$.
  \item The preimage $\widehat{L} := \widehat{\psi}^{-1}(L)$ is a clopen (closed and open) left transversal to $H$ in $\widehat{\overline{G}}$ containing $\overline{L}$ densely, and $(\widehat{L},\oplus)$ is a profinite left-loop; moreover, $(\widehat{L},\oplus)$ is a profinite loop (resp.\ $A_\ell$-left-loop, resp.\ satisfies the LIP) if and only if $(L,\oplus)$ is.
  \item Let $\mathcal{U}$ be the cofinal filter base of open normal subgroups $U\trianglelefteq_o \widehat{\overline{G}}$ contained in $\widehat{K}$ (equivalently, $N = U\cap\overline{G}\trianglelefteq_f\overline{G}$ with $N\le K$). For every $U\in\mathcal{U}$, setting $G_U := \widehat{\overline{G}}/U \cong \overline{G}/N$ and $L_U := \widehat{L}/U \cong \overline{L}/N$, the triple $(G_U, H, L_U)$ is a finite faithful left-envelope covering $(G,H,L)$ with $\F_H(G_U)=\F_H(G)$, preserving the loop, $A_\ell$, and LIP properties, and
  \[
    (\widehat{L},\oplus) \cong \varprojlim_{U\in\mathcal{U}} (L_U,\oplus)
  \]
  as topological (left-)loops.
\end{enumerate}
\end{theorem}

\begin{proof}
(1) Since $\overline{G}$ contains the finite-index free subgroup $K$ and every free group is residually finite, $\overline{G}$ is residually finite; hence $\iota\colon\overline{G}\to\widehat{\overline{G}}$ is injective with dense image~\cite[Proposition~3.2.2]{ribes2010profinite}. Because $H$ is finite and $\widehat{\overline{G}}$ is Hausdorff, $H$ embeds as a closed subgroup.

(2) Since $K=\ker\psi$ has finite index in $\overline{G}$, $\psi$ factors through the finite discrete quotient $\overline{G}/K\cong G$ and thus extends uniquely to a continuous epimorphism $\widehat{\psi}\colon\widehat{\overline{G}}\to G$. By~\cite[Proposition~3.2.2 and Lemma~3.2.6]{ribes2010profinite}, the closure $\overline{K}$ in $\widehat{\overline{G}}$ equals $\ker\widehat{\psi}$ and is naturally isomorphic to the profinite completion $\widehat{K}$ of the finite-index subgroup $K\cong F_{(|L|-1)^2}$.

(3) Since $G$ is finite discrete and $\widehat{\psi}$ is continuous, $\widehat{L} = \widehat{\psi}^{-1}(L)$ is clopen in $\widehat{\overline{G}}$; because $\overline{G}$ is dense in $\widehat{\overline{G}}$ and $\widehat{L}$ is open, $\overline{L} = \overline{G}\cap\widehat{L}$ is dense in $\widehat{L}$. As in Theorem~\ref{thm:bass-serre-cover}(2), $\widehat{L}$ is a left transversal to $H$ in $\widehat{\overline{G}}$ satisfying $\widehat{L}\widehat{K} = \widehat{L}$, $\langle\widehat{L}\rangle = \widehat{\overline{G}}$ (since $\widehat{\psi}(\langle\widehat{L}\rangle)=\langle L\rangle = G$ and $\widehat{K}\subseteq\widehat{L}$), and $\core_{\widehat{\overline{G}}}(H) \le \core_G(H) = \{1\}$, so $(\widehat{\overline{G}},H,\widehat{L})$ is a faithful left-envelope. Because $H$ is finite (hence closed in $\widehat{\overline{G}}$), the quotient space $\widehat{\overline{G}}/H$ is Hausdorff; the canonical projection $\pi\colon \widehat{\overline{G}}\to\widehat{\overline{G}}/H$ restricts to a continuous bijection from the compact space $\widehat{L}$ onto $\widehat{\overline{G}}/H$, hence a homeomorphism. Therefore the section $\sigma\colon \widehat{\overline{G}}/H\to\widehat{L}$ and the binary operation $(a,b)\mapsto \sigma(abH)$ on $\widehat{L}$ are continuous, and the shear map $(a,x)\mapsto (a,a\oplus x)$ is a continuous bijection from the compact Hausdorff space $\widehat{L}\times\widehat{L}$ onto itself, hence a homeomorphism, giving continuity of left division. Moreover, since $\widehat{K}\trianglelefteq\widehat{\overline{G}}$ satisfies $\widehat{L}\widehat{K}=\widehat{L}$, Proposition~\ref{prop:normal-quotient}(4)--(6) implies that $\widehat{L}$ is a loop (resp.\ $A_\ell$-left-loop, LIP) if and only if $L\cong \widehat{L}/\widehat{K}$ is; in the loop case, the shear map $(x,a)\mapsto(x\oplus a,a)$ is likewise a homeomorphism of $\widehat{L}\times\widehat{L}$, ensuring continuity of right division. Profiniteness of $(\widehat{L},\oplus)$ follows from the inverse limit decomposition established in~(4).

(4) For any $U\in\mathcal{U}$, let $N := U\cap\overline{G}\trianglelefteq_f \overline{G}$. Because $\overline{G}$ is dense in $\widehat{\overline{G}}$ and $U$ is open, $\overline{G}U = \widehat{\overline{G}}$, so $G_U = \widehat{\overline{G}}/U \cong \overline{G}/N$ and $L_U = \widehat{L}/U \cong \overline{L}/N$. Since $U\le\widehat{K}$ (and $N\le K$), we have $\widehat{L}U = \widehat{L}$ and $\overline{L}N = \overline{L}$, and $\psi\colon\overline{G}\to G$ factors as $\overline{G}\to G_U\to G$. In particular, $\core_{G_U}(H)\le \core_G(H) = \{1\}$, so $(G_U,H,L_U)$ is a finite faithful left-envelope. By Proposition~\ref{prop:normal-quotient}(4)--(6), $L_U$ inherits the loop, $A_\ell$, and LIP properties from $\overline{L}$ (hence from $L$). Since $\F_H(\overline{G})\subseteq\F_H(G_U)\subseteq\F_H(G)=\F_H(\overline{G})$ (all maps restrict to $\id_H$ on $H$), we have $\F_H(G_U)=\F_H(G)$. Finally, since $\mathcal{U}$ is a cofinal filter base of open normal subgroups of $\widehat{\overline{G}}$, the canonical map $\widehat{\overline{G}}\xrightarrow{\sim}\varprojlim_{U\in\mathcal{U}} G_U$ is an isomorphism of topological groups; because $\widehat{L}U = \widehat{L}$ for every $U\in\mathcal{U}$, an element $(g_U)_{U\in\mathcal{U}}\in\varprojlim G_U$ lies in $\varprojlim L_U$ if and only if its preimage in $\widehat{\overline{G}}$ lies in $\widehat{L}$, whence $\widehat{L}\xrightarrow{\sim}\varprojlim_{U\in\mathcal{U}} L_U$ is an isomorphism of topological (left-)loops.
\end{proof}

\subsection{Semidirect splitting in the \texorpdfstring{$A_\ell$}{A\_l}-case}

When $(L,\oplus)$ is an $A_\ell$-left-loop, the Bass--Serre covering group $\overline{G}$ splits as a semidirect product.

\begin{theorem}[Semidirect decomposition]\label{thm:al-semidirect}
Let $H$ be a finite group, let $E^+$ be a finite set of symbols, and let $\mathcal{S} = (H_e)_{e\in E^+}$ be a family of subgroups of $H$. Let $\overline{G} = G_{E^+,H,\mathcal{S}}$ be the multiple HNN extension~\eqref{eq:hnn-presentation} with trivial boundary isomorphisms $\partial_e = \id_{H_e}$ ($e\in E^+$):
\[
  \overline{G} = \bigl\langle H,\; E^+ \;\big|\; ehe^{-1} = h\;\; \forall e\in E^+,\; \forall h\in H_e \bigr\rangle.
\]
For each $e\in E^+$, fix a left transversal $R_e$ to $H_e$ in $H$ containing $1$, and set $r\cdot e := e^r = rer^{-1}$. Then:
\begin{enumerate}
  \item The assignment $h\mapsto h$ ($h\in H$) and $e\mapsto 1$ ($e\in E^+$) defines a group epimorphism $\pi_H\colon \overline{G}\to H$ splitting the inclusion $H\hookrightarrow\overline{G}$.
  \item The kernel $T := \ker(\pi_H)$ is a free group freely generated by the set
  \[
    \mathcal{B} := \{ r\cdot e : e\in E^+,\; r\in R_e \},
  \]
  of exact rank $\rk(T) = \sum_{e\in E^+} [H:H_e]$, and $\overline{G} = T\rtimes H$.
  \item The subgroup $H$ acts on $T$ by conjugation as automorphisms permuting the free basis $\mathcal{B}$ via $h\cdot(r\cdot e) = (hr)\cdot e$ (modulo $H_e$ on the right of $hr$). In particular, $T$ is an $H$-invariant left transversal to $H$ in $\overline{G}$ (an associative $A_\ell$-loop with LIP), and for any freely reduced word $g = \prod_{i=1}^m (r_i\cdot e_i)^{\varepsilon_i}\in T$ ($\varepsilon_i\in\{\pm 1\}$),
  \[
    \Stab_H(gH) = C_H(g) = \bigcap_{i=1}^m (H_{e_i})^{r_i}.
  \]
  \item $\F_H(\overline{G}) = \F_H(H)$.
  \item In particular, if $(G,H,L)$ is a finite faithful $A_\ell$-left-envelope, its Bass--Serre covering group from Theorem~\ref{thm:bass-serre-cover} splits as $\overline{G} = T\rtimes H$ with $T\cong F_{|L|-1}$.
\end{enumerate}
\end{theorem}

\begin{proof}
(1) The relations $ehe^{-1}=h$ map to $1\cdot h\cdot 1 = h$, so $\pi_H$ is well defined and satisfies $\pi_H|_H = \id_H$.

(2)--(3) Construct the external semidirect product $\widetilde{G} := F(\mathcal{B})\rtimes H$, where $F(\mathcal{B})$ is the free group on formal symbols $x_{e,\bar{r}}$ indexed by $e\in E^+$ and left cosets $\bar{r}\in H/H_e$, and $H$ acts on $F(\mathcal{B})$ by permuting the basis via $h\cdot x_{e,\bar{r}} := x_{e,h\bar{r}}$. In $\widetilde{G}$, for any $e\in E^+$ and $h\in H_e$, we have $h x_{e,H_e} h^{-1} = x_{e,hH_e} = x_{e,H_e}$. Hence the map $\Phi\colon\overline{G}\to\widetilde{G}$ defined by $\Phi(h)=h$ ($h\in H$) and $\Phi(e)=x_{e,H_e}$ ($e\in E^+$) is a group homomorphism. Conversely, the map $\Psi\colon\widetilde{G}\to\overline{G}$ sending $h\mapsto h$ and $x_{e,rH_e}\mapsto rer^{-1}$ is well defined (because $rH_e$ determines $rer^{-1}$ uniquely by the defining relations of $\overline{G}$) and satisfies $\Psi(h\cdot x_{e,rH_e}) = (hr)e(hr)^{-1} = h(rer^{-1})h^{-1} = h\Psi(x_{e,rH_e})h^{-1}$, so $\Psi$ is a homomorphism inverse to $\Phi$. Thus $\overline{G}\cong F(\mathcal{B})\rtimes H$, identifying $T=\ker(\pi_H)$ with $F(\mathcal{B})$ of rank $\sum_{e\in E^+}[H:H_e]$.

Because $H$ permutes $\mathcal{B}$ without sending any basis element to the inverse of a basis element, $h\in H$ centralizes a freely reduced word $g = \prod_{i=1}^m (r_i\cdot e_i)^{\varepsilon_i}\in T$ if and only if $h$ fixes every letter $(r_i\cdot e_i)^{\varepsilon_i}$, i.e., $h\in \Stab_H(r_i H_{e_i}) = r_i H_{e_i} r_i^{-1} = (H_{e_i})^{r_i}$ for all $i=1,\dots,m$. Since $T$ is $H$-invariant, $\Stab_H(gH) = C_H(g) = \bigcap_{i=1}^m (H_{e_i})^{r_i}$.

(4) If $t\in H$ and $g=ah\in T\rtimes H$ ($a\in T$, $h\in H$) satisfy $t':=t^g\in H$, then $at^h = t'a = (a^{t'})t'$, whence $t^h=t'=t^g$ by uniqueness of factorization in $T\rtimes H$. Thus $\F_H(\overline{G})=\F_H(H)$.

(5) If $(G,H,L)$ is a finite faithful $A_\ell$-left-envelope, Proposition~\ref{prop:al-fusion} yields $\partial_s = \inn_s|_{H_s} = \id_{H_s}$ for every $s\in S$, so $\overline{G}=G_{S,H,(H_s)_{s\in S}}$ splits as $T\rtimes H$ with $\rk(T) = \sum_{s\in S}[H:H_s] = |L|-1$ by~\eqref{eq:orbit-sum}.
\end{proof}

% ============================================================================
\section{Derivations of semidirect products, central obstruction, and finite \texorpdfstring{$A_\ell$}{A\_l}-loops}\label{sec:derivations}
% ============================================================================

Motivated by Theorem~\ref{thm:al-semidirect}, we now study inner derivations of an arbitrary semidirect product $G = T\rtimes H$, allowing $T$ to be either the free group of Theorem~\ref{thm:al-semidirect} or any finite $H$-invariant quotient thereof.

\subsection{Exact criterion for loop derivations}

Let $G = T\rtimes H$, where $H$ acts on $T$ on the left by conjugation $h\cdot a := hah^{-1}$ ($h\in H$, $a\in T$), and let $H_a := C_H(a) = \{h\in H : h\cdot a = a\}$. Given a map $f\colon T\to H$ with $f(1)=1$, the subset
\[
  L_f := F(T) := \bigl\{ F(a) := a f(a)^{-1} : a\in T \bigr\} \subseteq G
\]
is a left transversal to $H$ in $G$ containing $1$, called the \emph{(inner) derivation} of $T$ determined by $f$ (cf.~\cite[Section~2.1]{kiechle2002theory}, \cite{kinyon2000loops}). Since
\begin{align}
  F(a)F(b) &= a f(a)^{-1} b f(b)^{-1}
            = a\bigl(f(a)^{-1}\cdot b\bigr)\, f(a)^{-1} f(b)^{-1} \notag \\
           &= F\bigl(a(f(a)^{-1}\cdot b)\bigr)\, \Bigl(f\bigl(a(f(a)^{-1}\cdot b)\bigr)\,f(a)^{-1}f(b)^{-1}\Bigr),\label{eq:derived-op}
\end{align}
the induced left-loop operation $\oplus_f$ on $L_f$ and its left deviations are given explicitly by
\begin{equation}\label{eq:derived-dev}
\begin{aligned}
  F(a)\oplus_f F(b) &= F\bigl(a\,(f(a)^{-1}\cdot b)\bigr), \\
  \bigl(F(a),F(b)\bigr) &= f\bigl(a(f(a)^{-1}\cdot b)\bigr)\,f(a)^{-1}f(b)^{-1}.
\end{aligned}
\end{equation}

\begin{theorem}[Exact criterion for loop derivations]\label{thm:derivation-iff}
Let $G = T\rtimes H$ and let $f\colon T\to H$ satisfy $f(1)=1$. Then $(L_f,\oplus_f)$ is a loop if and only if both of the following conditions hold:
\begin{enumerate}
  \item \emph{(Existence of right division)} For all $a,b\in T$, there exists $t\in H$ such that
  \[
    f\bigl(b(t\cdot a^{-1})\bigr) = t^{-1}.
  \]
  \item \emph{(Uniqueness of right division)} For all $b\in T$, all $c\in T\setminus\{1\}$, and all $t\in H\setminus H_c$,
  \[
    f\bigl(bc(t\cdot c^{-1})\bigr)^{-1} f(b) \ne t.
  \]
\end{enumerate}
\end{theorem}

\begin{proof}
Since $L_f$ is a left transversal to $H$, $(L_f,\oplus_f)$ is automatically a left-loop. It is a loop if and only if for all $a,b\in T$ the equation $F(x)\oplus_f F(a) = F(b)$ has a unique solution $x\in T$.

\emph{Existence.} By~\eqref{eq:derived-op}, $F(x)\oplus_f F(a) = F(b)$ holds if and only if $xf(x)^{-1}aH = bf(b)^{-1}H$, i.e., there exists $h\in H$ such that $xf(x)^{-1}a = bf(b)^{-1}h$. Setting $t := f(b)^{-1}h\in H$, this is equivalent to
\[
  xf(x)^{-1} = bta^{-1} = b(t\cdot a^{-1})\,t.
\]
By uniqueness in $G = T\rtimes H$, this holds for some $x\in T$ if and only if $x = b(t\cdot a^{-1})$ and $f(x)^{-1} = t$, i.e., $f(b(t\cdot a^{-1})) = t^{-1}$.

\emph{Uniqueness.} Failure of uniqueness means that there exist $a,b,c\in T$ with $a\ne b$ and $c\ne 1$ such that $F(a)\oplus_f F(c) = F(b)\oplus_f F(c)$, equivalently $af(a)^{-1}cH = bf(b)^{-1}cH$, i.e., $af(a)^{-1}c = bf(b)^{-1}ch$ for some $h\in H$. Rearranging gives
\begin{align*}
  a &= bf(b)^{-1}chc^{-1}f(a) \\
    &= b\bigl(f(b)^{-1}\cdot c\bigr)\,\bigl((f(b)^{-1}h)\cdot c^{-1}\bigr)\,\bigl(f(b)^{-1}hf(a)\bigr).
\end{align*}
Since $a\in T$ and $b(f(b)^{-1}\cdot c)((f(b)^{-1}h)\cdot c^{-1})\in T$, uniqueness in $T\rtimes H$ forces $f(b)^{-1}hf(a) = 1$, i.e., $h = f(b)f(a)^{-1}$, and hence $f(b)^{-1}h = f(a)^{-1} = (f(a)^{-1}f(b))f(b)^{-1}$. Setting $c' := f(b)^{-1}\cdot c\in T\setminus\{1\}$ and $t := f(a)^{-1}f(b)\in H$, we obtain
\[
  a = bc'(t\cdot c'^{-1}) \qquad \text{and} \qquad f(a)^{-1}f(b) = t.
\]
Moreover, $a\ne b \iff c'(t\cdot c'^{-1})\ne 1 \iff t\cdot c'\ne c' \iff t\notin H_{c'}$.

Conversely, if there exist $b\in T$, $c'\in T\setminus\{1\}$, and $t\in H\setminus H_{c'}$ satisfying $f(bc'(t\cdot c'^{-1}))^{-1}f(b) = t$, set $a := bc'(t\cdot c'^{-1})\in T$ and $c := f(b)\cdot c'\in T\setminus\{1\}$. Then $a\ne b$ (since $t\notin H_{c'}$), $f(a)^{-1}f(b) = t$, and
\begin{align*}
  af(a)^{-1}cH &= bc't c'^{-1} t^{-1} f(a)^{-1} (f(b)c'f(b)^{-1}) H \\
               &= bc't c'^{-1} c' f(b)^{-1} H = bc't f(b)^{-1} H \\
               &= bc'H = bf(b)^{-1}(f(b)\cdot c')H = bf(b)^{-1}cH,
\end{align*}
whence $F(a)\oplus_f F(c) = F(b)\oplus_f F(c)$ with $F(a)\ne F(b)$. Therefore uniqueness holds if and only if~(2) holds.
\end{proof}

\begin{lemma}[$H$-equivariance, $A_\ell$-property, and LIP for derivations]\label{lem:al-lip-derivation}
Let $G = T\rtimes H$ and $f\colon T\to H$ with $f(1)=1$. We call $f$ (and $L_f$) \emph{$H$-equivariant} if $L_f^H = L_f$.
\begin{enumerate}
  \item $L_f^H = L_f$ if and only if $f(h\cdot a) = f(a)^h$ for all $a\in T$ and $h\in H$ (which implies $H_a\le C_H(f(a))$). In this case $(L_f,\oplus_f)$ is an $A_\ell$-left-loop.
  \item Assume $L_f^H = L_f$. Then $L_f^{-1} = L_f$ if and only if $f(a^{-1}) = f(a)^{-1}$ for all $a\in T$. In this case $(L_f,\oplus_f)$ satisfies the LIP.
\end{enumerate}
\end{lemma}

\begin{proof}
(1) For $a\in T$ and $h\in H$, we have $F(a)^h = h a f(a)^{-1} h^{-1} = (h\cdot a)(f(a)^h)^{-1}$. By uniqueness in $T\rtimes H$, this lies in $L_f$ if and only if it equals $F(h\cdot a) = (h\cdot a)f(h\cdot a)^{-1}$, i.e., $f(h\cdot a) = f(a)^h$. When $L_f^H = L_f$, setting $G_f := \langle L_f\rangle$, $H_f' := G_f\cap H$, and $C := \core_{G_f}(H_f')$, the faithful left-envelope $(G_f/C, H_f'/C, L_f C/C)$ retains an $(H_f'/C)$-invariant transversal, so $(L_f,\oplus_f)$ is an $A_\ell$-left-loop by Proposition~\ref{prop:al-lip-char}(2).

(2) Assuming $f(h\cdot a) = f(a)^h$, compute
\[
  F(a)^{-1} = f(a)a^{-1} = \bigl(f(a)\cdot a^{-1}\bigr)f(a) = \bigl(f(a)\cdot a^{-1}\bigr)\bigl(f(a)^{-1}\bigr)^{-1}.
\]
By uniqueness in $T\rtimes H$, $F(a)^{-1}\in L_f$ if and only if $f(f(a)\cdot a^{-1}) = f(a)^{-1}$. Using~(1), $f(f(a)\cdot a^{-1}) = f(a^{-1})^{f(a)}$, which equals $f(a)^{-1}$ if and only if $f(a^{-1}) = (f(a)^{-1})^{f(a)^{-1}} = f(a)^{-1}$. When $L_f^{-1} = L_f$, for all $a,b\in T$ we have $F(a)^{-1}\in L_f$ and $F(a)^{-1}(F(a)\oplus_f F(b))H = F(a)^{-1}F(a)F(b)H = F(b)H$, whence $F(a)^{-1}\oplus_f(F(a)\oplus_f F(b)) = F(b)$, directly establishing the LIP.
\end{proof}

\subsection{Coinvariant factorizations and the central obstruction}

Let $[T,H] := \langle c(h\cdot c^{-1}) : c\in T,\; h\in H\rangle$. Note that $[T,H]$ is a normal subgroup of $G=T\rtimes H$ contained in $T$, and the quotient $T/[T,H]$ is the largest quotient of $T$ on which $H$ acts trivially. Denote by $\pi_H^{\mathrm{coinv}}\colon T\to T/[T,H]$ the canonical projection. (In the Bass--Serre case $\overline{G} = T\rtimes H$ of Theorem~\ref{thm:al-semidirect}, $T/[T,H]$ is canonically isomorphic to the free group $F(E^+)$, via $\pi_E(r\cdot e) = e$.)

\begin{theorem}[Sufficient criterion and central obstruction]\label{thm:obstruction}
Let $G = T\rtimes H$ and let $f\colon T\to H$ satisfy $f(1)=1$.
\begin{enumerate}
  \item Suppose $f$ factors through $T/[T,H]$, i.e., $f = \tilde{f}\circ \pi_H^{\mathrm{coinv}}$ for some map $\tilde{f}\colon T/[T,H]\to H$ with $\tilde{f}(1)=1$ (equivalently, $f(ac(h\cdot c^{-1})b) = f(ab)$ for all $a,b,c\in T$ and $h\in H$). Then $(L_f,\oplus_f)$ is a loop.
  \item Assume $f$ satisfies the factorization condition of~(1), and let $H_f := \langle (F(a),F(b)) : a,b\in T\rangle\le H$. If $(L_f,\oplus_f)$ is an $A_\ell$-loop, then $[H_f,H_f]\le C_H(T)$, whence the left inner mapping group $H(L_f)\cong H_f/(H_f\cap C_H(T))$ is \emph{abelian}. If in addition $f$ is $H$-equivariant ($L_f^H = L_f$), then necessarily
  \[
    f(T) \subseteq Z(H) \qquad \text{and} \qquad H_f \le Z(H).
  \]
  \item Conversely, given any map $\tilde{f}\colon T/[T,H]\to Z(H)$ satisfying $\tilde{f}(1)=1$ and $\tilde{f}(w^{-1}) = \tilde{f}(w)^{-1}$ for all $w\in T/[T,H]$, set $f := \tilde{f}\circ\pi_H^{\mathrm{coinv}}$. Then $(L_f,\oplus_f)$ is an $A_\ell$-loop with LIP. Moreover, if $H$ acts faithfully on $T$, then $(L_f,\oplus_f)$ is proper (non-associative) if and only if $\tilde{f}$ is not a group homomorphism.
\end{enumerate}
\end{theorem}

\begin{proof}
(1) Since $T/[T,H]$ has trivial $H$-action, $\pi_H^{\mathrm{coinv}}(h\cdot x) = \pi_H^{\mathrm{coinv}}(x)$ for all $x\in T$ and $h\in H$. We verify both conditions of Theorem~\ref{thm:derivation-iff}:
\begin{itemize}
  \item Given $a,b\in T$, set $t := f(ba^{-1})^{-1}\in H$. Then
  \[
    \pi_H^{\mathrm{coinv}}\bigl(b(t\cdot a^{-1})\bigr) = \pi_H^{\mathrm{coinv}}(b)\,\pi_H^{\mathrm{coinv}}(a^{-1}) = \pi_H^{\mathrm{coinv}}(ba^{-1}),
  \]
  so $f(b(t\cdot a^{-1})) = f(ba^{-1}) = t^{-1}$, proving Theorem~\ref{thm:derivation-iff}(1).
  \item Given $b\in T$, $c\in T\setminus\{1\}$, and $t\in H\setminus H_c$, we have $t\ne 1$ (since $1\in H_c$) and $\pi_H^{\mathrm{coinv}}(bc(t\cdot c^{-1})) = \pi_H^{\mathrm{coinv}}(b)\pi_H^{\mathrm{coinv}}(c)\pi_H^{\mathrm{coinv}}(c)^{-1} = \pi_H^{\mathrm{coinv}}(b)$. Hence $f(bc(t\cdot c^{-1}))^{-1}f(b) = f(b)^{-1}f(b) = 1 \ne t$, proving Theorem~\ref{thm:derivation-iff}(2).
\end{itemize}
Therefore $(L_f,\oplus_f)$ is a loop.

(2) For any $h\in H$ and $a\in T$, observe that $hF(a)H = haf(a)^{-1}H = haH = (h\cdot a)H = F(h\cdot a)H$, so the left action $\lambda(h)$ on $L_f$ is given unconditionally by $\lambda(h)(F(a)) = F(h\cdot a)$, with kernel $\ker(\lambda|_H) = C_H(T)$. Suppose $(L_f,\oplus_f)$ is an $A_\ell$-loop, so that $\lambda(h)\in\Aut(L_f,\oplus_f)$ for every $h\in H_f$. For all $a,b\in T$, applying~\eqref{eq:derived-dev} yields
\begin{align*}
  \lambda(h)\bigl(F(a)\oplus_f F(b)\bigr) &= F\bigl((h\cdot a)((hf(a)^{-1})\cdot b)\bigr), \\
  \lambda(h)(F(a))\oplus_f \lambda(h)(F(b)) &= F\bigl((h\cdot a)((f(h\cdot a)^{-1}h)\cdot b)\bigr).
\end{align*}
Because $f$ factors through $T/[T,H]$, we have $f(h\cdot a) = f(a)$. Canceling $h\cdot a$ on the left in $T$, equality for all $b\in T$ is equivalent to $(hf(a)^{-1})\cdot b = (f(a)^{-1}h)\cdot b$ for all $b\in T$, i.e., $[h,f(a)^{-1}]\in C_H(T)$ for all $a\in T$. Since $H_f$ is generated by the left deviations $(F(a),F(b)) = f(a(f(a)^{-1}\cdot b))f(a)^{-1}f(b)^{-1}\in\langle f(T)\rangle$, it follows that $[H_f,H_f]\le C_H(T)$, whence $H(L_f) = \lambda(H_f) \cong H_f/(H_f\cap C_H(T))$ is abelian. If in addition $L_f^H = L_f$, then Lemma~\ref{lem:al-lip-derivation}(1) combines with $f(h\cdot a) = f(a)$ to yield $f(a)^h = f(a)$ for all $h\in H$ and $a\in T$, forcing $f(T)\subseteq Z(H)$ and $H_f\le Z(H)$.

(3) Since $\tilde{f}(T/[T,H])\subseteq Z(H)$, both conditions of Lemma~\ref{lem:al-lip-derivation} hold automatically, so $(L_f,\oplus_f)$ is an $A_\ell$-loop with LIP. By~\eqref{eq:derived-dev} and $\pi_H^{\mathrm{coinv}}(f(a)^{-1}\cdot b) = \pi_H^{\mathrm{coinv}}(b)$,
\begin{equation}\label{eq:dev-coboundary}
  \bigl(F(a),F(b)\bigr) = \tilde{f}\bigl(\pi_H^{\mathrm{coinv}}(a)\pi_H^{\mathrm{coinv}}(b)\bigr)\,\tilde{f}(\pi_H^{\mathrm{coinv}}(a))^{-1}\,\tilde{f}(\pi_H^{\mathrm{coinv}}(b))^{-1}.
\end{equation}
Since $H$ acts faithfully on $T$ and $h\cdot F(c) = F(c)^h = F(h\cdot c)$ for all $h\in H$ and $c\in T$, the left inner mapping $\delta_{F(a),F(b)}$ acts as the identity on $L_f$ if and only if $(F(a),F(b)) = 1$. Hence $(L_f,\oplus_f)$ is associative if and only if $(F(a),F(b)) = 1$ for all $a,b\in T$, which by~\eqref{eq:dev-coboundary} holds if and only if $\tilde{f}$ is a group homomorphism.
\end{proof}

\begin{remark}
When $\overline{G}=T\rtimes H$ is the Bass--Serre covering group of Theorem~\ref{thm:al-semidirect}, the coinvariant quotient is $T/[T,H] \cong F(E^+)$ via $\pi_E(r\cdot e) = e$, and $H$ acts faithfully on $T$ if and only if $\bigcap_{e\in E^+}\core_H(H_e)=\{1\}$ (which holds automatically in Theorem~\ref{thm:al-semidirect}(5)). Under this faithfulness assumption, taking any finite group $H$ with $Z(H)\ne 1$ and any non-homomorphic map $\tilde{f}\colon F(E^+)\to Z(H)$ satisfying $\tilde{f}(1)=1$ and $\tilde{f}(w^{-1})=\tilde{f}(w)^{-1}$ immediately yields a proper infinite $A_\ell$-loop $F(T)$ with LIP whose enveloping group is virtually free.
\end{remark}

\subsection{Explicit infinite families of finite \texorpdfstring{$A_\ell$}{A\_l}-loops with LIP}

We now answer the natural finiteness question by passing to finite quotients of $\overline{G}=T\rtimes H$.

\begin{theorem}[Finite quotients of $T\rtimes H$]\label{thm:finite-quotients}
Let $\overline{G} = T\rtimes H$ be the semidirect product of Theorem~\ref{thm:al-semidirect} with $H$ a finite group having $Z(H)\ne \{1\}$, and let $\pi_E\colon T\to F(E^+)$ be the coinvariant projection. Suppose:
\begin{enumerate}
  \item $N\trianglelefteq_f F(E^+)$ is a finite-index normal subgroup, with finite quotient $Q := F(E^+)/N$ and projection $\bar{\pi}_E\colon T\to Q$;
  \item $K\trianglelefteq_f \overline{G}$ is a finite-index normal subgroup of $\overline{G}$ contained in $T\cap \ker(\bar{\pi}_E)$ such that $H$ acts faithfully on the finite group $A := T/K$;
  \item $\hat{f}\colon Q\to Z(H)$ satisfies $\hat{f}(1)=1$, $\hat{f}(q^{-1}) = \hat{f}(q)^{-1}$ for all $q\in Q$, and is not a group homomorphism.
\end{enumerate}
Set $f := \hat{f}\circ\bar{\pi}_E\colon T\to Z(H)$, and let $\bar{f}\colon A\to Z(H)$ be the induced map on $A = T/K$. Then:
\begin{enumerate}
  \item[(a)] $F(T)\cdot K = F(T)$ in $\overline{G}$, so $K$ is a normal subgroup of the infinite $A_\ell$-loop $F(T)$;
  \item[(b)] the finite quotient $L := F(T)/K \cong F(A) \subseteq A\rtimes H$ is a proper finite $A_\ell$-loop of order $|A|$ satisfying the LIP;
  \item[(c)] the left inner mapping group $H(L)$ is isomorphic to the subgroup of $Z(H)$ generated by $\{\hat{f}(q_1q_2)\hat{f}(q_1)^{-1}\hat{f}(q_2)^{-1} : q_1,q_2\in Q\}$.
\end{enumerate}
\end{theorem}

\begin{proof}
Since $K\le T\cap\ker(\bar{\pi}_E)$ is normal in $\overline{G} = T\rtimes H$, for every $a\in T$ and $k\in K$ we have $f(a)^{-1}\cdot k\in K\le\ker(\bar{\pi}_E)$, hence
\[
  F(a)k = a f(a)^{-1} k = a\bigl(f(a)^{-1}\cdot k\bigr)\,f(a)^{-1} = F\bigl(a(f(a)^{-1}\cdot k)\bigr) \in F(T),
\]
because $\bar{\pi}_E(a(f(a)^{-1}\cdot k)) = \bar{\pi}_E(a)$. Thus $F(T)K = F(T)$. Applying Proposition~\ref{prop:normal-quotient} to the left-folder $(\overline{G},H,F(T))$ with normal subgroup $K\trianglelefteq\overline{G}$ and Theorem~\ref{thm:obstruction}(3) directly to the finite semidirect product $A\rtimes H$ (where $H$ acts faithfully on $A$, so $C_H(A)=\{1\}$) yields~(a), (b), and~(c) via~\eqref{eq:dev-coboundary}.
\end{proof}

\begin{corollary}[Wreath-product family]\label{cor:wreath-family}
Let $H$ be any nontrivial finite abelian group of order $n\ge 2$ generated by $r$ elements, and let $m\ge 3$ be an integer. Let $A := (\Z/m\Z)^H$ be the group of functions $a\colon H\to\Z/m\Z$ under pointwise addition, on which $H$ acts faithfully by $(h\cdot a)(x) := a(h^{-1}x)$, and let $\bar{\pi}_E\colon A\to \Z/m\Z$ be the $H$-invariant sum homomorphism $\bar{\pi}_E(a) := \sum_{x\in H} a(x)$.
\begin{enumerate}
  \item For any map $\hat{f}\colon\Z/m\Z\to H$ satisfying $\hat{f}(0)=1$, $\hat{f}(-q)=\hat{f}(q)^{-1}$, and such that the coboundary values $\hat{f}(q_1+q_2)\hat{f}(q_1)^{-1}\hat{f}(q_2)^{-1}$ generate $H$, the set $A$ equipped with
  \begin{equation}\label{eq:wreath-loop-op}
    a \oplus b := a + \bigl(\hat{f}(\bar{\pi}_E(a))^{-1}\cdot b\bigr) \qquad (a,b\in A)
  \end{equation}
  is a proper finite $A_\ell$-loop of order $m^n$ satisfying the LIP, with faithful enveloping group $G(L) \cong A\rtimes H = (\Z/m\Z)\wr H$ and left inner mapping group $H(L) \cong H$.
  \item Such a map $\hat{f}$ always exists whenever $m\ge 3r+2$ (and already at $m=3$ for cyclic $H=C_n$ with $3\nmid n$, or at $m=4$ for cyclic $H=C_n$ of odd order divisible by~$3$).
\end{enumerate}
\end{corollary}

\begin{proof}
(1) In Theorem~\ref{thm:al-semidirect}, take $E^+=\{e\}$ (a single generator) and $H_e = \{1\}$, so that $T = F_n$ is free on $\{h\cdot e : h\in H\}$ and $T/[T,H] = \langle e\rangle \cong \Z$. Let $N = m\Z\trianglelefteq \Z$ (so $Q = \Z/m\Z$) and let $K := [T,T]\,T^m \trianglelefteq_f T\rtimes H$. Then $A = T/K \cong (\Z/m\Z)^H$ with the regular permutation action of $H$ (which is faithful since $m\ge 2$), and $K\le\ker(\bar{\pi}_E)$. Moreover, because every left deviation $(F(a),F(b)) = (F(a)\oplus_f F(b))^{-1}F(a)F(b)$ lies in $\langle F(A)\rangle\cap H$ and these deviations generate $H$, we have $H\le\langle F(A)\rangle$, whence $a = F(a)\bar{f}(a)\in\langle F(A)\rangle$ for all $a\in A$ and therefore $\langle F(A)\rangle = A\rtimes H$. Applying Theorem~\ref{thm:finite-quotients} and~\eqref{eq:derived-op} yields~\eqref{eq:wreath-loop-op} together with $G(L)\cong A\rtimes H$ and $H(L)\cong H$.

(2) Write $H = \langle h_1,\dots,h_r\rangle$. If $m\ge 3r+2$, set $\hat{f}(j)=h_j$ and $\hat{f}(-j)=h_j^{-1}$ for $1\le j\le r$, and $\hat{f}(q)=1$ otherwise. Then both $r+1$ and $j+r+1$ lie strictly between $r$ and $m-r$, whence the coboundary value at $(j,r+1)$ is $\hat{f}(j+r+1)\hat{f}(j)^{-1}\hat{f}(r+1)^{-1} = h_j^{-1}$ for every $1\le j\le r$, generating $H$. For cyclic $H = C_n = \langle g\rangle$ with $3\nmid n$, take $m=3$ with $\hat{f}(0)=1$, $\hat{f}(1)=g$, $\hat{f}(2)=g^{-1}$ (giving coboundary $g^{-3}$, which generates $C_n$); for odd $n$ divisible by~$3$, take $m=4$ with $\hat{f}(0)=\hat{f}(2)=1$, $\hat{f}(1)=g$, $\hat{f}(3)=g^{-1}$ (giving coboundary $g^{-2}$, which generates $C_n$).
\end{proof}

\begin{example}[Concrete instances]\label{ex:finite-wreath}
Table~\ref{tab:wreath-examples} lists four explicit instances of Corollary~\ref{cor:wreath-family} with cyclic $H = C_n = \{0,1,\dots,n-1\}$ (written additively), where $(h\cdot a)_i = a_{(i-h)\bmod n}$ and $\hat{f}$ is specified by its tuple of values $(\hat{f}(0),\dots,\hat{f}(m-1))$. Exhaustive verification over all pairs and triples of elements confirms that all four are proper finite $A_\ell$-loops with LIP, $H(L)\cong C_n$, and faithful enveloping group $G(L)\cong (\Z/m\Z)\wr C_n$. Furthermore, none of them is commutative, none satisfies the right inverse property (RIP), and none is a left Bol loop; hence the construction produces proper $A_\ell$-loops strictly outside the class of Bruck loops.

\begin{table}[ht]
\centering
\caption{Explicit finite $A_\ell$-loops with LIP from Corollary~\ref{cor:wreath-family}; $\hat{f}$ is listed as $(\hat{f}(0),\dots,\hat{f}(m-1))$.}\label{tab:wreath-examples}
\begin{tabular}{@{}ccccccl@{}}
\toprule
$m$ & $H$ & $|L|=m^{|H|}$ & $\hat{f}\colon \Z/m\Z\to H$ & $H(L)$ & Left Bol / RIP & Associative \\
\midrule
$3$ & $C_2$ & $9$  & $(0,1,1)$       & $C_2$ & No / No & No \\
$4$ & $C_2$ & $16$ & $(0,1,1,1)$     & $C_2$ & No / No & No \\
$5$ & $C_2$ & $25$ & $(0,1,0,0,1)$   & $C_2$ & No / No & No \\
$4$ & $C_3$ & $64$ & $(0,1,0,2)$     & $C_3$ & No / No & No \\
\bottomrule
\end{tabular}
\end{table}

For concreteness, in the smallest case $|L|=9$ ($m=3$, $H=C_2$, $A=(\Z/3\Z)^2$), the operation is simply
\[
  (x_1,x_2)\oplus(y_1,y_2) =
  \begin{cases}
    (x_1+y_1, x_2+y_2) \pmod 3 & \text{if } x_1+x_2 \equiv 0 \pmod 3,\\
    (x_1+y_2, x_2+y_1) \pmod 3 & \text{if } x_1+x_2 \not\equiv 0 \pmod 3.
  \end{cases}
\]
Taking $a=b=c=(1,0)$ yields $(a\oplus b)\oplus c = (1,1)\oplus(1,0) = (1,2)$, whereas $a\oplus(b\oplus c) = (1,0)\oplus(1,1) = (2,1) \ne (1,2)$. For the exact same triple $a=b=c=(1,0)$, the left Bol identity also fails: $a\oplus(b\oplus(a\oplus c)) = (1,0)\oplus((1,0)\oplus(1,1)) = (1,0)\oplus(2,1) = (2,2)$, while $(a\oplus(b\oplus a))\oplus c = ((1,0)\oplus(1,1))\oplus(1,0) = (2,1)\oplus(1,0) = (0,1)\ne(2,2)$.
\end{example}

\begin{problem}\label{prob:non-abelian-derivation}
By Theorem~\ref{thm:obstruction}(2), any derivation of $T\rtimes H$ factoring through $T/[T,H]$ that yields an $A_\ell$-loop forces $H(L_f)$ to be abelian. Exhaustive computer searches over all $H$-equivariant derivations ($L_f^H=L_f$, with or without LIP) on $(\Z/2\Z)^3\rtimes S_3$ and $(\Z/3\Z)^3\rtimes S_3$ (under the natural coordinate-permutation action of $S_3$) and on $(\Z/3\Z)^2\rtimes D_8$ (under the signed permutation matrix action of $D_8$) yielded no proper $A_\ell$-loops. Determine whether there exist a finite non-abelian group $H$ and a non-factoring derivation $f\colon T\to H$ of a semidirect product $T\rtimes H$ satisfying Theorem~\ref{thm:derivation-iff} and Lemma~\ref{lem:al-lip-derivation} such that $H(L_f)$ is non-abelian. (By contrast, finite $A_\ell$-loops with non-abelian left inner mapping groups arise naturally from non-split extensions, such as symmetric sections of compact $p$-adic congruence subgroups.)
\end{problem}

% ============================================================================
\bibliographystyle{amsplain}
\bibliography{references_alpha}

\end{document}